\documentclass[10pt,a4paper]{article} 
\usepackage{latexsym,amsfonts,theorem} 
\usepackage{graphicx}

\theoremstyle{break}
\theorembodyfont{\slshape}
\newtheorem{Thm}{Theorem}[section]

\newtheorem{Lem}[Thm]{Lemma}
\newtheorem{Cor}[Thm]{Corollary}

\theorembodyfont{\upshape}
\newtheorem{Def}[Thm]{Definition}
\newtheorem{Rem}[Thm]{Remark}

\newcommand{\BoP}[1]{\noindent {\sc Proof#1: }}
\newcommand{\EoP}{\hfill$\Box$\vspace{6pt}}
\newcommand{\NoP}{\hfill$\Box$}
\newcommand{\bb}{\mathbb}

\newcommand{\reins}{{\bb R}}
\newcommand{\R}{{\cal R}}
\newcommand{\N}{{\bb N}}

\title{Fourier knots}

\author{\\ Christoph Lamm  \\  \\ 
{\footnotesize R\"{u}ckertstr.~3,} \\ 
{\footnotesize 65187 Wiesbaden, Germany,} \\ 
{\footnotesize e-mail: christoph.lamm@web.de}}

\begin{document}
\parindent0.5cm
\date{}

\maketitle

\vspace{1.5cm}
\begin{quotation}
\abstract{We show that every knot has a checkerbord diagram and that
every knot is the closure of a rosette braid. 
We define Fourier knots of type $(n_1,n_2,n_3)$ as knots which have
parametrizations where each coordinate function $x_i (t)$ is a finite
Fourier series of length $n_i$, and conclude that
every knot is a Fourier knot of type $(1,1,n)$ for some natural
number $n$.
}
\end{quotation}

\section{Rosette braids}
Let $B_n$ be the braid group on $n$ strings, $\pi: B_n \to S_n$ be the
map to the symmetric group on $n$ letters and $P_n =$ ker$(\pi)$ be the 
pure braid group on $n$ strings. The generators of $P_n$ are denoted
by $A_{i,j}$. If $\pi_0$ is a permutation, then $k(\pi_0)$ denotes its number
of cycles.

\begin{Def} 
A braid of the form
$$
\prod_{i=1}^n \Bigl[\prod_{j \,{\rm odd} \atop j < s} 
\sigma_j^{\varepsilon_{i,j}}
\prod_{j \,{\rm even} \atop j < s} 
\sigma_j^{\varepsilon_{i,j}}\Bigr]\,,\,\,\varepsilon_{i,j} \in \{\pm 1\}
$$
is called a \emph{rosette braid of type (s,n)}.
The set of rosette braids of type $(s,n)$ is denoted by $\R (s,n)$.
\end{Def}

\begin{Lem}\label{l1}
\begin{enumerate}
\item[a)]
\begin{itemize}
\item[i)]
$\alpha \in \R(s,1) \Rightarrow k(\pi(\alpha))=1.$
\item[ii)]
$\alpha \in \R(s,s) \Rightarrow \pi(\alpha)$ =id.
\item[iii)]
$\alpha \in \R(s,ns+1) \Rightarrow k(\pi(\alpha))=1.$
\end{itemize}
\item[b)]
For each generator $A_{i,j}$ of $P_n$ there is an
$\alpha \in \R(s,s)$, so that $A_{i,j} = \alpha$.
\end{enumerate}
\end{Lem}

\BoP{}
a) Proposition i)  is true for $s=2$, because then the braid
word has the form $\sigma^{\pm 1}$. An element of $\R (s,1)$
is built out of an element of $\R (s-1,1)$ by a Markov-II-move
(insertion of $\sigma_{s-1}^{\pm 1}$). Because the number of
components is unchanged by a Markov-II-move, we conclude by induction
that the proposition holds for all $s$.

ii) As shown in part i), the permutation $\pi(\alpha)$ of a braid
$\alpha \in \R(s,1)$ consists of one cycle. Hence the permutation
of a braid in $\R(s,s)$ is the trivial permutation on $s$ letters.
Part iii) is an immediate consequence.

b) We consider the two strings $i$ and $j$ ($i < j$) of the 
braid $A_{i,j}$. If as above
$\pi_1 = \pi(\alpha)$ is the permutation of a braid $\alpha \in \R(s,1)$,
then $\pi_1$ is a cycle of length $s$. Hence there is a $k$ with $1 \le k
< s$, so that $\pi_1^k(i) > \pi_1^k(j)$, and thus the strings $i$ and $j$ 
cross each other. It is possible to arrange the strings in such a way that all
strings but $i$ and $j$ can be pulled tight and the strings $i,j$ form
the generators $A_{i,j}$ or $A_{i,j}^{-1}$.
\EoP

\begin{Thm}
Let $\alpha \in B_s$ be a braid, with closure a knot.
Then $\alpha$ is conjugate to a rosette braid of type
$(s,ns+1)$ for a suitable $n$.
\end{Thm}

\BoP{}
Let $\alpha \in B_s$ be a braid, so that $\hat \alpha$ is
a knot. Let $\delta$ be an arbitrary braid in $\R(s,1)$ and
$\pi_1$ its permuation. 
The permutations $\pi(\alpha)$ and $\pi_1$ are conjugate in the
symmetric group because both consist of one cycle. 
Let $\beta \in B_s$ be a braid, so that
$\pi(\beta)^{-1}\pi(\alpha)\pi(\beta)=\pi_1$. 
Then $\delta^{-1}\beta^{-1}\alpha\beta$ is a pure braid and
because of Lemma \ref{l1} we can write it as an element of
$\R(s,ns)$ for a suitable $n$.
Multiplication with $\delta$ yields $\beta^{-1}\alpha\beta$
as an element of $\R(s,ns+1)$. Hence we have shown that 
$\alpha$ is conjugate to a rosette braid of type $(s,ns+1)$.
\EoP

The braid index of a knot $K$ is denoted by $br(K)$.

\begin{Cor}
Every knot $K$ is the closure of a rosette braid with $br(K)$ strings.
\NoP
\end{Cor}

\bigskip

\begin{figure}
\centerline{\includegraphics[width=0.5\textwidth]{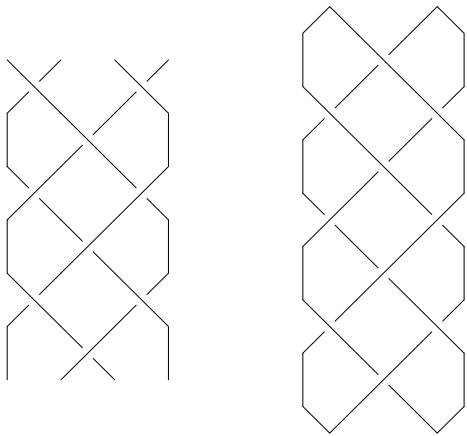}}
\caption{\small A rosette braid of type $(4,3)$ and a checkerboard diagram of type $(4,3)$}
\label{r1}
\end{figure}


\section{Checkerboard diagrams}
\begin{Def}
A knot diagram is called a \emph{checkerboard diagram of type
$(2b,n)$}, if it is the plat closure of a braid
$\sigma_2^{\varepsilon_2}\ldots\sigma_{2b-2}^{\varepsilon_{2b-2}}\cdot \alpha$
with $\alpha \in \R(2b,n)$ and $\varepsilon_2,\ldots,\varepsilon_{2b-2}\in 
\{\pm 1\}$.
\end{Def}

\begin{figure}
\centerline{\includegraphics[width=0.7\textwidth]{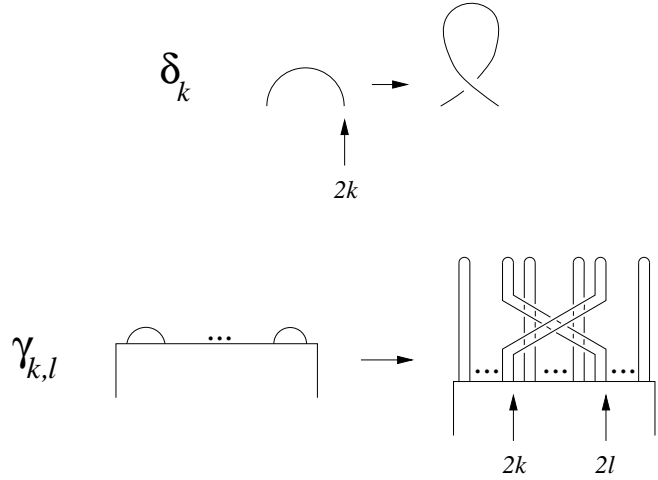}}
\caption{\small $\delta$ and $\gamma$ to modify 
the permutation of a braid which constitutes a plat. The $\delta'$
and $\gamma'$ are the mirrored operations for the lower plat closure.}
\label{op}
\end{figure}

Let $\pi_0$ be the permutation of the braid 
$\sigma_2^{\varepsilon_2}\ldots\sigma_{2b-2}^{\varepsilon_{2b-2}}\in B_{2b}$.
The plat-operations $\delta, \gamma, \delta$ and $\gamma'$ which we need
for Lemma \ref{oper} are defined in Figure \ref{op}.

\bigskip

\begin{Lem}\label{oper}
Let $K$ be a knot which is given as a plat closure of a braid
$\alpha \in B_s$. Then there is a sequence of operations 
$\delta$, $\delta'$, $\gamma$ and $\gamma'$ which transforms the plat
$\overline \alpha$ to a plat $\overline \beta$ with $\pi (\beta)= \pi_0$. 
\end{Lem}

\BoP{}
Using the operations of Figure \ref{op}, the permutation $\pi_0$ can
be produced step by step. We start with string 1 and move its end-position
to the position $\pi_0(1)=1$. Then, travelling along the knot we
can successively adjust the end-positions at the upper and lower
plat-closure of the braid. The result is a plat with permutation
$\pi_0$.
\EoP

\bigskip

\begin{figure}
\centerline{\includegraphics[width=0.7\textwidth]{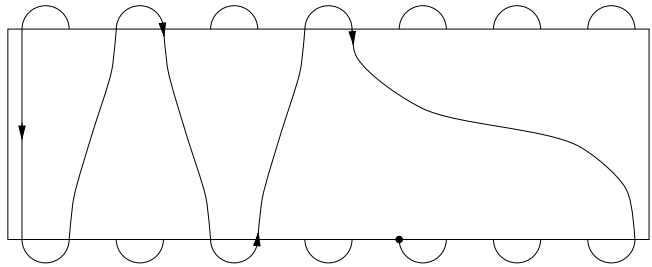}}
\caption{\small Illustration for the proof of Lemma \ref{oper}. The operation
$\gamma_{5,7}$ followed by $\delta_5$ adjust the string of the plat which is
travelled as the fifth string.}
\label{op2}
\end{figure}

\begin{Thm}\label{ch}
Every knot with bridge number $b$ has a checkerboard diagram of
type $(2b,2nb)$ for a suitable $n$.
\end{Thm}

\BoP{}
In Lemma \ref{oper} we succeeded to represent the knot $K$ as a plat
$\bar \alpha$ with $\pi(\alpha)=\pi_0$.
We consider the pure braid 
$\beta=\sigma_2\sigma_4\ldots\sigma_{s-2}\cdot\alpha$. 
By Lemma \ref{l1} the braid $\beta$ can be written as a rosette 
braid $\beta'$ of type $(s,ns)$ for some natural number $n$. 
Hence $\sigma_2^{-1}\sigma_4^{-1}\ldots\sigma_{s-2}^{-1}\cdot \beta'$
is a checkerboard diagram for $K$. If we choose the plat representative
of $K$ with $2b(K)$ strings, then the checkerboard diagram is of type
$(2b(K),2nb(K))$.
\EoP

\section{Fourier knots}

We call a series of the form
$$
\sum_{i=1}^n \alpha_i \cos(2\pi m_i t + \varphi_i)
$$
with $t, \alpha_i, \varphi_i \in \reins, m_i \in \N \,(i=1,\ldots,n)$
a \emph{finite Fourier series of length $n$}.

\begin{Def}
A knot is a \emph{Fourier knot of type $(n_1,n_2,n_3)$} if it can be parametrized
by coordinate functions $x_1, x_2, x_3 : [0,1]\to \reins$ which are finite Fourier series  
of length $n_1\le n_2 \le n_3$.
\end{Def}

\begin{Rem}
The Fourier knots of type $(1,1,1)$ are the \emph{Lissajous knots}. They 
were studied in the articles \cite{Bogle}, \cite{JP} and \cite{Lamm}.
Not all knots are Lissajous knots, but the next theorem proclaims that
every knot is a Fourier knot of an especially simple type. Fourier knots
were also defined in \cite{Kau} and \cite{Trau}. 
\end{Rem}

\begin{Thm}
Every knot $K$ is a Fourier knot of type $(1,1,n_K)$ for some $n_K\in {\bb N}$.
\end{Thm}

\BoP{}
We consider a Fourier knot of type $(1,1,n)$ and its projection on the
$x$-$y$-plane. By \cite{Lamm} the knot diagram is a checkerboard diagram.
Conversely, by Theorem \ref{ch} every knot has a checkerboard diagram. 
The height-function in $z$-direction can be approximated by a finite
Fourier series.
\EoP

\begin{Rem}
The trefoil knots and the figure-eight knot are not Lissajous knots. 
They are Fourier knots of type $(1,1,2)$. The parametrizations of 
\cite{Bogle} (with a correction of misprints) are
 
\centerline{
$x_1(t)=\cos(2t+6)$, $x_2(t)=\cos(3t+0.15)$, $x_3(t)=\cos(4t+1)+\cos(5t)$} 

\noindent
for a trefoil and

\centerline{
$x_1(t)=\cos(2t+0.8)$, $x_2(t)=\cos(3t+0.15)$, $x_3(t)=\cos(4t+1)+\cos(5t)$}

\noindent
for the figure-eight knot. Here we use the interval $t\in[0,2\pi]$, in
order to have the same parametrization as in \cite{Bogle}.
\end{Rem}

\begin{Def}
If $K$ is a knot, the \emph{Fourier index of $K$} is 
the smallest number $n$ for which $K$ is a Fourier
knot of type $(1,1,n)$. 
\end{Def}

It is not known if there are knots with arbitrarily high Fourier index.

\bigskip
{\bf Acknowledgement:} 
I thank J.\,Kneissler for discussions on braids and plats.

\end{document}